\documentclass[11pt,a4paper]{article}
\usepackage[T1]{fontenc}
\usepackage{lmodern}
\usepackage[margin=28mm]{geometry}
\usepackage{amsmath,amssymb,amsthm,mathtools}
\usepackage{enumitem,booktabs,microtype,xcolor,needspace}
\usepackage[colorlinks=true,linkcolor=blue!50!black,citecolor=blue!50!black,urlcolor=blue!50!black]{hyperref}
\hypersetup{pdftitle={Normalized solutions of critical biharmonic equations on bounded domains},pdfauthor={Qilin Xie and Yelin Chen}}
\numberwithin{equation}{section}
\newtheorem{theorem}{Theorem}[section]
\newtheorem{proposition}[theorem]{Proposition}
\newtheorem{lemma}[theorem]{Lemma}
\newtheorem{corollary}[theorem]{Corollary}
\theoremstyle{remark}
\newtheorem{remark}[theorem]{Remark}
\newcommand{\R}{\mathbb R}
\newcommand{\dd}{\,dx}
\newcommand{\norm}[1]{\left\|#1\right\|}
\newcommand{\weak}{\rightharpoonup}
\newcommand{\dist}{\operatorname{dist}}
\newcommand{\cl}{\operatorname{cl}}
\newcommand{\eps}{\varepsilon}
\newcommand{\loc}{\mathrm{loc}}
\allowdisplaybreaks[2]
\setlist[enumerate]{label=\textup{(\roman*)},leftmargin=*}
\title{Normalized solutions of critical biharmonic equations\\ on bounded domains}
\author{Qilin Xie\thanks{Corresponding author. This work was supported by the Guangdong Basic and Applied Basic Research Foundation (Grant No. 2026A1515012458).} \quad Yelin Chen\\[0.5em]
\small Guangdong University of Technology,\\
\small School of Mathematics and Statistics,\\
\small Guangzhou, Guangdong 510520, P.R. China}
\date{}
\begin{document}
\maketitle
\begin{abstract}
We study the critical biharmonic equation with clamped boundary conditions
and prescribed $L^2$ norm on a smooth bounded domain in $\R^N$, $N\ge5$.
For every sufficiently small mass, the constrained energy admits a local
minimizer. We also obtain a second solution for masses in the closure of
a measurable set having positive measure in every interval adjacent to
zero. This holds on arbitrary smooth bounded domains when $N\ge8$ and,
when $N=5,6,7$, under a condition involving the first clamped eigenvalue
and a cutoff quotient. Explicit sufficient conditions are given and
verified on balls. The proof combines a quantitative monotonicity
argument with a norm estimate excluding loss of compactness.
\end{abstract}
\noindent\textbf{Keywords.} Biharmonic equation; normalized solution;
critical Sobolev exponent; constrained mountain pass.

\medskip
\noindent\textbf{Mathematics Subject Classification (2020):}
Primary 35J91; Secondary 35J35, 35J30.

\section{Introduction}
Let $\Omega\subset\R^N$, $N\ge5$, be a bounded connected domain with
smooth boundary. We study the following equation,
\begin{equation}\label{eq:original}
\begin{cases}
\Delta^2U=\lambda U+|U|^{q-2}U&\text{in }\Omega,\\
U=\partial_\nu U=0&\text{on }\partial\Omega,\\
\displaystyle\int_\Omega U^2\dd=\rho^2,
\end{cases}
\end{equation}
where $q=\frac{2N}{N-4}$ and both the real-valued functions $U$ and $\lambda\in\R$ are unknown. All solutions are weak solutions in $H_0^2(\Omega)$, with clamped boundary
conditions understood in the trace sense.
The operator $\Delta^2$ is the principal part of the linear plate
operator. In the clamped model, the displacement and its normal
slope vanish on the boundary.  For background on these issues, including eigenvalues, Green
functions and semilinear critical problems, we refer to Gazzola,
Grunau and Sweers \cite{GGSw}.

The critical exponent is the endpoint of the continuous embedding
$H_0^2(\Omega)\hookrightarrow L^q(\Omega)$. Concentration prevents
compactness at this endpoint, even on a bounded domain. The classical
work of Br\'ezis and Nirenberg \cite{BN} explains how a lower-order
perturbation interacts with critical growth in second-order equations.
The fourth-order analogue raises related questions, but with different
scaling and a more delicate dependence on the boundary conditions.
Critical biharmonic equations on bounded domains have been studied in
connection with the fourth-order analogue of the Br\'ezis--Nirenberg
problem. In the clamped case, the corresponding prescribed-parameter
problem is
\begin{equation}\label{eq:fixedparameter}
\begin{cases}
\Delta^2v=\lambda v+|v|^{q-2}v&\text{in }\Omega,\\
v=\partial_\nu v=0&\text{on }\partial\Omega,
\end{cases}
\end{equation}
where $\lambda$ is fixed. Edmunds, Fortunato and Jannelli \cite{EFJ}
proved the existence of a nontrivial solution for
$0<\lambda<\Lambda_1$ when $N\ge8$. For $N=5,6,7$, they obtained
existence on an interval immediately below $\Lambda_1$, with refined
radial results on balls. Their analysis relates the dimension threshold
to the sharp Sobolev constant and to estimates for truncated
whole-space solutions. Gu, Deng and Wang \cite{GDW} studied existence
and nonexistence for the same critical clamped problem. Using a
Palais--Smale condition below the critical Sobolev energy and a
mountain-pass argument, they obtained dimension-dependent conditions
on the prescribed parameter, including results on balls in dimensions
five to seven.

The geometry of the domain and the choice of boundary conditions also
play a role. Gazzola, Grunau and Squassina \cite{GGS} considered critical
growth problems with both clamped and Navier boundary conditions. Their
results include existence on certain contractible domains obtained by
small-capacity perturbations of domains with nontrivial topology, as
well as compactness and nonexistence results. Singular lower-order
terms lead to a further dimensional dependence. Xiong and Shen
\cite{XS} studied clamped problems in which the linear term is
$\lambda|x|^{-s}v$, $0\le s\le4$, with either subcritical or critical
power growth. They obtained existence results and investigated the
change in critical dimensions as the strength of the singularity varies.

For Navier boundary conditions, Cheng, Feng and Wei \cite{CFW}
considered biharmonic equations with weighted concave-convex
nonlinearities on smooth bounded domains. In the subcritical range,
they used Nehari-manifold decompositions, fibering maps and minimization
to obtain one or two nontrivial solutions under conditions on the
weights. The associated space is $H^2(\Omega)\cap H_0^1(\Omega)$,
whereas the clamped problem is posed in $H_0^2(\Omega)$.
Thus the boundary conditions in these two variational settings must be
distinguished, even though both energies involve $\norm{\Delta v}_2^2$.

Throughout, $\rho$ is the
$L^2$ norm; its square is the integral mass. Define
\begin{equation}\label{eq:E}
E(U)=\frac12\norm{\Delta U}_2^2-\frac1q\norm U_q^q,
\qquad S_\rho=\{U\in H_0^2(\Omega):\norm U_2=\rho\}.
\end{equation}
Let
\begin{equation}\label{eq:constants}
S_2=\inf_{0\ne v\in D^{2,2}(\R^N)}
\frac{\norm{\Delta v}_{L^2(\R^N)}^2}{\norm v_{L^q(\R^N)}^2},\qquad
\Lambda_1=\inf_{0\ne v\in H_0^2(\Omega)}
\frac{\norm{\Delta v}_2^2}{\norm v_2^2}>0,
\end{equation}
where $D^{2,2}(\R^N)$ is the completion of $C_c^\infty(\R^N)$ under
$\norm{\Delta\cdot}_2$.

In equation \eqref{eq:original}, the $L^2$ norm is prescribed and $\lambda$ arises
as a Lagrange multiplier. Its value and sign are therefore part of the
problem. Prescribed-parameter existence results such as those for
equation \eqref{eq:fixedparameter} do not directly determine the masses of the
solutions. Moreover, critical points must now be constructed on an
$L^2$ sphere, on which the critical energy is unbounded below.
The loss of compactness at the critical exponent persists under this
constraint.

The unboundedness below can be seen directly from an interior
concentration. Choose $\phi\in C_c^\infty(B_1)$ with
$\norm\phi_{L^2(\R^N)}=1$, fix an interior point $x_0$, and set
$U_t(x)=\rho t^{N/2}\phi(t(x-x_0))$ for sufficiently large $t$.
Then $U_t\in S_\rho$ and
\[
E(U_t)=\frac{\rho^2t^4}{2}\norm{\Delta\phi}_{L^2(\R^N)}^2
-\frac{\rho^qt^{4N/(N-4)}}{q}\norm\phi_{L^q(\R^N)}^q
\longrightarrow-\infty.
\]
A first solution must therefore be sought as a local minimum rather
than a global minimum on $S_\rho$. For small $\rho$, the estimate
\[
E(U)\ge\frac12\norm{\Delta U}_2^2
-\frac{S_2^{-q/2}}q\norm{\Delta U}_2^q
\]
provides a positive barrier around a low-norm region. Beyond this
barrier, concentration supplies the descending portion of a
mountain-pass path. The difficulty is to turn this geometry into a
critical point at the critical exponent.

The normalized formulation has its own variational history. Noris,
Tavares and Verzini \cite{NTV} studied prescribed-mass ground states
and orbital stability for Sobolev-subcritical NLS on bounded domains.
Their work highlights the role of local minimization when the energy
is unbounded below on the mass constraint. They subsequently studied
normalized systems, including Sobolev-critical nonlinearities
\cite{NTVsystems}. This places the first solution in a broader
variational framework; the construction of a second solution requires
additional control of constrained mountain-pass sequences.

Normalized fourth-order problems have also been studied in the whole
space. Liu, Sun and Zhang \cite{LSZ} consider mixed dispersion and
combined nonlinearities involving Sobolev-critical growth. Their
analysis distinguishes the construction of radial critical sequences
from recovery of the prescribed mass. On a bounded domain, compactness
in $L^2$ preserves the mass of a bounded weakly convergent sequence,
but concentration in the critical norm remains possible. Thus the
bounded-domain setting removes one obstacle while retaining the
essential difficulty associated with critical growth.
For further related work on biharmonic Schr\"odinger equations,
we refer to \cite{TSW,WY,LTW,LY}.

For second-order equations on bounded domains, Pierotti, Verzini and
Yu \cite{PVY} developed a quantitative monotonicity argument for the
Sobolev-critical normalized problem on general domains. Song and Zou
\cite{SZ} established positive normalized ground-state and mountain-pass
solutions on star-shaped bounded domains. These works motivate the present study
of the critical clamped biharmonic equation. We use the
parameter-derivative method of \cite{PVY}; compactness is recovered
through a strict bound on the derivative norm of the constrained
Palais--Smale sequence.

The monotonicity approach is rooted in Jeanjean's construction of
bounded Palais--Smale sequences \cite{J}. Here its quantitative form
is needed because boundedness alone does not exclude concentration.
The parameter derivative must provide a strict norm bound compatible
with the fourth-order Sobolev threshold. The corresponding splitting
uses the Br\'ezis--Lieb lemma \cite{BL}; the sharp Sobolev constant
and fourth-order optimizers are discussed in \cite{Sw}.
The fourth-order problem requires two further comparisons. First,
a nonzero concentration defect costs at least $S_2^{N/4}$ in squared
Laplacian norm after scaling back to the original equation. The weak
limit retains the prescribed $L^2$ mass and costs at least
$\Lambda_1\rho^2$. We obtain an upper bound strictly below the sum
of these two quantities, which excludes concentration without using
the sign of a Pohozaev boundary term. Second, the $L^2$ norm of a
truncated Sobolev optimizer changes its leading order at $N=8$.
In the lower dimensions the error surviving normalization must be
compared with the first eigenvalue. This leads to an additional
cutoff condition, which can be checked explicitly on balls and on
some non-star-shaped domains.

The result below separates the two existence assertions accordingly.
The local minimum is obtained at every sufficiently small mass on
any smooth bounded domain. The second solution is obtained on the
closure of a selected measurable set of masses accumulating at zero.
For $N\ge8$ no additional geometric condition is needed; for
$N=5,6,7$ the sufficient condition is stated next. The estimates and
examples following the theorem describe this distinction explicitly.

For $5\le N\le7$, we use the following additional hypothesis: for some $x_0\in\Omega$ there exists
$\eta\in C_c^\infty(\Omega)$, $0\le\eta\le1$, equal to one near $x_0$,
such that, setting $f(x)=|x-x_0|^{4-N}$ and extending $\eta$ by zero,
\begingroup
\renewcommand{\theHequation}{hypothesis.H}
\begin{equation}\label{eq:cutoffhyp}
\begin{gathered}
D_\eta=\int_{\R^N}|\Delta((1-\eta)f)|^2\dd,\qquad
L_\eta=\int_\Omega\eta^2f^2\dd,\\
\frac{D_\eta}{2L_\eta}<\frac{2\Lambda_1}{N}.
\end{gathered}\tag{H}
\end{equation}
\endgroup
The integrals are finite and positive. The cutoff is chosen once,
independently of $\rho$ or $\mu$.
For $N\ge8$, this hypothesis is not required.

\Needspace{9\baselineskip}
\begin{theorem}\label{thm:main}
Let $N\ge5$ and let $\Omega\subset\R^N$ be a smooth bounded connected
domain. The following assertions hold.
\begin{enumerate}
\item \textup{First solution.} There exists $\rho_*>0$ such that, for
every $0<\rho<\rho_*$, equation \eqref{eq:original} has a solution
$(U_\rho^{\loc},\lambda_\rho^{\loc})$. The function $U_\rho^{\loc}$
is a local minimizer of $E|_{S_\rho}$ and satisfies
\begin{equation}\label{eq:localmain}
\frac2N\Lambda_1\rho^2\le E(U_\rho^{\loc})
<\frac{\Lambda_1\rho^2}{2}<\frac2N S_2^{N/4}.
\end{equation}
\item \textup{Second solution.} Assume either $N\ge8$, or
$N\in\{5,6,7\}$ and \eqref{eq:cutoffhyp}. In the latter case, fix
$C_0$ such that
\begin{equation}\label{eq:C0choice}
\frac{D_\eta}{2L_\eta}<C_0<\frac{2\Lambda_1}{N}.
\end{equation}
After possibly decreasing $\rho_*$, there exist
$\sigma\in(0,\Lambda_1)$ and a measurable set $Q\subset(0,\rho_*)$ with
\begin{equation}\label{eq:Qmeasure}
|Q\cap(0,t)|>0\qquad\text{for every }0<t<\rho_*.
\end{equation}
For each $\rho\in\overline Q\cap(0,\rho_*)$, equation
\eqref{eq:original} has a solution $(U_\rho^2,\lambda_\rho^2)$ with
$U_\rho^2\ne\pm U_\rho^{\loc}$, and
\begin{equation}\label{eq:secondmain}
\frac2N S_2^{N/4}\le E(U_\rho^2)\le\frac2N S_2^{N/4}+
\begin{cases}
C\rho^{(N-4)/2},&N>8,\\
C\rho^2/|\log\rho|,&N=8,\\
C_0\rho^2,&N=5,6,7,
\end{cases}
\end{equation}
\begin{equation}\label{eq:secondnorm}
\norm{\Delta U_\rho^2}_2^2
\le S_2^{N/4}+(\Lambda_1-\sigma)\rho^2.
\end{equation}
Here $C$ is independent of $\rho$, and one may choose
$\sigma=\Lambda_1/2$ for $N\ge8$. If $\rho\in Q$, or if
$\rho_j\uparrow\rho$ for some sequence $(\rho_j)\subset Q$, the second
solution may be chosen at the mountain-pass level $c_\rho$ defined in
Section~\ref{sec:geometry}.
\end{enumerate}
\end{theorem}

\begin{corollary}\label{cor:geometry}
For $N=5,6,7$, let $R_\Omega$ be the inradius of $\Omega$. Then the following results holds,
\begin{enumerate}
\item  Hypothesis
\eqref{eq:cutoffhyp} holds if
\begin{equation}\label{eq:domaincondition}
\Lambda_1(\Omega)R_\Omega^4>T_N,\qquad
T_5=\frac{3628625}{7334},\quad T_6=480,\quad T_7=\frac{4851}{19}.
\end{equation}
\item  Every ball satisfies this strict condition. Thus the full theorem
applies on every ball in every dimension $N\ge5$.
\item   Condition \eqref{eq:domaincondition} does not require star-shapedness.
\end{enumerate}

\end{corollary}

\begin{remark}\label{rem:scope}
Theorem~\ref{thm:main}\textup{(i)} requires no geometric assumption
other than boundedness and smoothness. For the second solution, no
star-shapedness assumption is needed when $N\ge8$. In dimensions
$5,6,7$, geometry enters only through the quantitative comparison
\eqref{eq:cutoffhyp}. In particular, the proof does not use a sign
condition on $x\cdot\nu$ or on the Lagrange multiplier.
This is relevant because on a general domain the boundary term in a
Pohozaev identity need not have a prescribed sign, and the usual
argument for bounding constrained Palais--Smale sequences by that
identity is then unavailable.

The low-dimensional conditions in Corollary~\ref{cor:geometry} can also
hold on domains with holes. Indeed, if
$B_R(x_0)\subset\Omega\subset B_{R'}(x_0)$, domain monotonicity and
scaling of the clamped eigenvalue give
\[
\Lambda_1(\Omega)R_\Omega^4
\ge\Lambda_1(B_{R'})R_\Omega^4
=\Lambda_1(B_1)(R_\Omega/R')^4
\ge\Lambda_1(B_1)(R/R')^4.
\]
Since $\Lambda_1(B_1)>T_N$, this is greater than $T_N$ when $R/R'$
is sufficiently close to one. A small interior hole can be placed in
$B_{R'}(x_0)\setminus\overline{B_R(x_0)}$ without changing these
inclusions. The resulting smooth connected domain is not star-shaped
with respect to any point. Subsection~\ref{sec:explicit} gives explicit
bounds for the permitted ratio $R'/R$.
Thus the condition is compatible with nontrivial topology; the
constants $T_N$ are sufficient thresholds obtained from specific
cutoffs, not necessary conditions for existence.

The two solutions have different variational origins. The first is
obtained for every sufficiently small mass. For the second, parameter
selection yields \eqref{eq:Qmeasure}, which guarantees arbitrarily
small admissible masses but does not imply that $Q$ has full measure
or is dense in $(0,\rho_*)$. Compactness extends the existence result
to $\overline Q\cap(0,\rho_*)$; identification with the mountain-pass
level is retained at the masses stated in Theorem~\ref{thm:main}.
Neither positivity nor a least-energy characterization among all
normalized solutions is asserted. In particular, the clamped space
is not in general preserved by $u\mapsto|u|$, so the second-order
absolute-value argument does not apply.
\end{remark}

\begin{remark}\label{rem:compactnessroute}
The critical embedding
$H_0^2(\Omega)\hookrightarrow L^{2N/(N-4)}(\Omega)$ is not compact.
In addition, a bound on the energy of a constrained Palais--Smale
sequence does not by itself give a bound for its $H_0^2$ norm.
We first write
\[
U=\rho u,\qquad \norm u_2=1,\qquad
\mu=\rho^{8/(N-4)},\qquad E(\rho u)=\rho^2J_\mu(u),
\]
so that the parameter appears in the functional on a fixed sphere.
The quantitative monotonicity construction gives, at differentiability
points of the mountain-pass level $d$,
\[
J_\mu(u_n)\to d(\mu),\qquad
(J_\mu|_M)'(u_n)\to0,\qquad
\limsup_n\norm{\Delta u_n}_2^2
\le2d(\mu)-2\mu d'(\mu)+6\mu.
\]
Here $M$ and $J_\mu$ are defined in Section~\ref{sec:prelim}.
The derivative term is controlled on a selected positive-measure set
by comparing $d$ with the explicit upper bound $g$ in
Lemma~\ref{lem:upper}. Since inequalities between functions cannot
be differentiated in general, Lemma~\ref{lem:selection} obtains
\[
-d'(\mu)\le(1+\delta(\mu))[-g'(\mu)]
\]
through an integral comparison over $(t,\gamma t)$.

Compactness then follows from a norm comparison. Indeed, a bounded
constrained Palais--Smale sequence has a weak limit $u$ with
$\norm u_2=1$, by the compact embedding into $L^2$. For $z_n=u_n-u$,
norm splitting and the limiting equation give
\[
\norm{\Delta z_n}_2^2-\mu\norm{z_n}_q^q\to0.
\]
If a subsequence has $\norm{\Delta z_n}_2^2\to A>0$, the sharp
Sobolev inequality implies
\[
A\ge S_2(A/\mu)^{2/q}
\quad\Longrightarrow\quad A\ge S_2^{N/4}\mu^{-(N-4)/4}.
\]
Since $\norm{\Delta u}_2^2\ge\Lambda_1$, loss of compactness would force
\[
\liminf_n\norm{\Delta u_n}_2^2
\ge S_2^{N/4}\mu^{-(N-4)/4}+\Lambda_1.
\]
The selected derivative estimate instead yields
\[
\limsup_n\norm{\Delta u_n}_2^2
\le S_2^{N/4}\mu^{-(N-4)/4}+\Lambda_1-\sigma,
\qquad\sigma>0,
\]
which excludes a nonzero defect. This argument uses neither a sign
condition on the multiplier nor a classification of concentration
profiles at interior or boundary points.
\end{remark}

\begin{remark}\label{rem:dimensionroute}
The additional assumption in dimensions $5,6,7$ comes from the
$L^2$ normalization of the test bubble. The whole-space Sobolev
optimizer has decay $W(x)\asymp|x|^{4-N}$, and hence
\[
\int_{1<|x|<R}W(x)^2\dd\asymp\int_1^R r^{7-N}\,dr.
\]
This integral is bounded for $N>8$, grows logarithmically for $N=8$,
and grows like $R^{8-N}$ for $5\le N\le7$. Consequently, for the
cutoff bubbles of Lemma~\ref{lem:bubble},
\[
L_\eps=\norm{\eta W_\eps}_2^2\asymp
\begin{cases}
\eps^4,&N>8,\\
\eps^4|\log\eps|,&N=8,\\
\eps^{N-4},&N=5,6,7.
\end{cases}
\]
The scales contributing to the maximum of the test path satisfy
$L_\eps\asymp\rho^2$. After division of the physical energy by
$\rho^2$, the cutoff error therefore tends to zero for $N\ge8$.
For $N=5,6,7$, it instead has a nonzero limit:
\[
\frac{A_\eps-S_2^{N/4}}{2L_\eps}
\longrightarrow\frac{D_\eta}{2L_\eta},\qquad
 d(\mu)\le\frac2N S_2^{N/4}\mu^{-(N-4)/4}
 +\frac{D_\eta}{2L_\eta}+o(1).
\]
Thus taking a smaller mass alone does not remove this error.

There are two restrictions on the constant used in parameter selection.
The upper energy estimate requires $C_0>D_\eta/(2L_\eta)$.
With $\delta(\mu)=d_0\mu^{(N-4)/4}$, the integral comparison in
Lemma~\ref{lem:selection} requires
\[
\frac{2d_0}{N}S_2^{N/4}>C_0,
\]
whereas exclusion of a concentration defect requires
\[
2C_0+\frac{N-4}{N}S_2^{N/4}d_0<\Lambda_1.
\]
These two strict inequalities admit a common $d_0$ exactly when
\[
\frac{NC_0}{2S_2^{N/4}}
<\frac{N(\Lambda_1-2C_0)}{(N-4)S_2^{N/4}}
\quad\Longleftrightarrow\quad C_0<\frac{2\Lambda_1}{N}.
\]
A permissible pair $(C_0,d_0)$ therefore exists precisely when the
chosen cutoff satisfies (H). This explains both the coefficient
$2/N$ in that hypothesis and the need for a strict inequality.

Subsection~\ref{sec:explicit} turns this analytic requirement into a
condition involving only the inradius and the first eigenvalue.
For the explicit limiting profiles, $ND_G/(4L_G)=T_N$, while scaling
an interior ball of radius $R$ gives
\[
\frac{D_{\eta_R}}{2L_{\eta_R}}
=R^{-4}\frac{D_\eta}{2L_\eta}.
\]
The strict inequality $\Lambda_1(\Omega)R_\Omega^4>T_N$ leaves enough
room to smooth the profiles without losing (H). Improving the cutoff
quotient would improve these sufficient thresholds; failure of a
particular threshold does not rule out (H), or the existence of a
second solution.
\end{remark}

Section~\ref{sec:prelim} contains the compactness criterion, and
Section~\ref{sec:geometry} establishes local minimization and the
mountain-pass geometry. Section~\ref{sec:construction} combines the
quantitative monotonicity argument, the upper energy bound and
parameter selection to obtain compact mountain-pass sequences.
Section~\ref{sec:results} proves Theorem~\ref{thm:main} and
Corollary~\ref{cor:geometry}, including the explicit verification of
the low-dimensional hypothesis.

\section{Preliminaries}\label{sec:prelim}
Equip $X=H_0^2(\Omega)$ with
$\langle u,v\rangle_X=\int_\Omega\Delta u\Delta v\dd$ and
$\norm u_X=\norm{\Delta u}_2$, an equivalent Hilbert norm.
The embedding $X\hookrightarrow L^2(\Omega)$ is compact. Extension by zero
maps $X$ into $D^{2,2}(\R^N)$, and therefore
\begin{equation}\label{eq:sobolev}
S_2\norm u_q^2 \le\norm{\Delta u}_2^2 \qquad \text{for\ any}\ u\in X.
\end{equation}
The sharp constant and the corresponding whole-space optimizers are
given in \cite{Sw}.

Set
\begin{equation}\label{eq:mu}
U=\rho u,\quad \mu=\rho^{q-2},\quad M=\{u\in X:\norm u_2=1\},
\quad J_\mu(u)=\frac12\norm{\Delta u}_2^2-\frac\mu q\norm u_q^q.
\end{equation}
The conversion formulas are
\begin{equation}\label{eq:conversion}
q-2=\frac8{N-4},\qquad \rho=\mu^{(N-4)/8},\qquad
\rho^2=\mu^{(N-4)/4},\qquad E(\rho u)=\rho^2J_\mu(u).
\end{equation}
The set $M$ is a closed smooth Hilbert submanifold of $X$, with
\begin{equation}\label{eq:tangent}
T_uM=\left\{\varphi\in X:\int_\Omega u\varphi\dd=0\right\}.
\end{equation}
The induced Riemannian metric is complete
and local manifold charts identify the intrinsic and ambient topologies
near that limit. Derivative norms on $T_uM$ are induced by $X$.

\begin{lemma}[Multipliers for bounded constrained sequences]\label{lem:multiplier}
Fix $\mu>0$. Suppose $(u_n)\subset M$ is bounded in $X$ and
$\norm{(J_\mu|_M)'(u_n)}\to0$. Set
\begin{equation}\label{eq:lambda}
\lambda_n=\norm{\Delta u_n}_2^2-\mu\norm{u_n}_q^q.
\end{equation}
Then $(\lambda_n)$ is bounded and
\begin{equation}\label{eq:residual}
\Delta^2u_n-\lambda_nu_n-\mu|u_n|^{q-2}u_n\longrightarrow0
\quad\text{in }X^*.
\end{equation}
Every weak limit $u$ of a subsequence belongs to $M$ and, after extracting
a further subsequence with $\lambda_n\to\lambda$, satisfies
\begin{equation}\label{eq:normalized}
\int_\Omega\Delta u\Delta\varphi\dd
=\lambda\int_\Omega u\varphi\dd+
\mu\int_\Omega|u|^{q-2}u\varphi\dd\qquad(\varphi\in X).
\end{equation}
\end{lemma}
\begin{proof}
For $\varphi\in X$, put
$P_n\varphi=\varphi-(\int_\Omega u_n\varphi\dd)u_n$.
Then $P_n\varphi\in T_{u_n}M$, and boundedness of $u_n$ and the continuous
embedding into $L^2$ give $\norm{P_n\varphi}_X\le C\norm\varphi_X$.
Since $J_\mu'(u_n)[u_n]=\lambda_n$,
\[
J_\mu'(u_n)[\varphi]-\lambda_n\int_\Omega u_n\varphi\dd
=J_\mu'(u_n)[P_n\varphi]=o_n(1)\norm\varphi_X.
\]
This proves \eqref{eq:residual}. Formula \eqref{eq:lambda} and the Sobolev
embedding give boundedness of $\lambda_n$.

Take a subsequence with $u_n\weak u$ in $X$, $u_n\to u$ in $L^2$ and
almost everywhere, and $\lambda_n\to\lambda$. Strong $L^2$ convergence
gives $\norm u_2=1$. The functions $|u_n|^{q-2}u_n$ are bounded in
$L^{q/(q-1)}$ and converge almost everywhere to $|u|^{q-2}u$, hence
converge weakly in that space. This assertion follows by uniform
integrability against bounded test functions, followed by density of
bounded functions in $L^q$ and the uniform dual norm bound.
Passing to the limit in \eqref{eq:residual} proves \eqref{eq:normalized}.
\end{proof}

\begin{lemma}\label{lem:compact}
Let $(u_n)\subset M$ be a bounded constrained Palais--Smale sequence for
$J_\mu$. Up to a subsequence, $u_n\weak u\in M$, and $u$ solves
\eqref{eq:normalized}. If this subsequence does not converge strongly,
a further subsequence satisfies
\begin{align}
\liminf_n\norm{\Delta u_n}_2^2
&\ge\norm{\Delta u}_2^2+S_2^{N/4}\mu^{-(N-4)/4}\ge\Lambda_1+S_2^{N/4}\mu^{-(N-4)/4},
\label{eq:normloss}\\
\liminf_n J_\mu(u_n)&\ge J_\mu(u)+{\frac2N S_2^{N/4}}\mu^{-(N-4)/4}.
\label{eq:energyloss}
\end{align}
In particular, such a sequence is precompact if
\begin{equation}\label{eq:normcriterion}
\limsup_n\norm{\Delta u_n}_2^2<S_2^{N/4}\mu^{-(N-4)/4}+\Lambda_1.
\end{equation}
\end{lemma}
\begin{proof}
Set $z_n=u_n-u$ and Hilbert space
orthogonality and the Brezis--Lieb lemma \cite{BL} give
\begin{align}
\norm{\Delta u_n}_2^2&=\norm{\Delta u}_2^2+\norm{\Delta z_n}_2^2+o_n(1),
\label{eq:hilbertsplit}\\
\norm{u_n}_q^q&=\norm u_q^q+\norm{z_n}_q^q+o_n(1).
\label{eq:BLsplit}
\end{align}
Testing the limiting equation with $u$, using \eqref{eq:lambda} and
$\norm u_2=1$, yields
\begin{equation}\label{eq:defectbalance}
\norm{\Delta z_n}_2^2-\mu\norm{z_n}_q^q\longrightarrow0.
\end{equation}
If convergence is not strong, choose a further subsequence with
$A=\lim\norm{\Delta z_n}_2^2>0$. Then \eqref{eq:sobolev} implies
\[
A\ge S_2(A/\mu)^{2/q}\qquad \text{and}\qquad
A\ge S_2^{N/4}\mu^{-(N-4)/4}.
\]
Since $u\in M$, $\norm{\Delta u}_2^2\ge\Lambda_1$, it follows that
\eqref{eq:normloss}. The energy splitting is
\[
J_\mu(u_n)=J_\mu(u)+\left(\frac12-\frac1q\right)A+o_n(1),
\]
which proves \eqref{eq:energyloss}. Under \eqref{eq:normcriterion},
every subsequence has a weakly convergent further subsequence and any
positive defect is excluded by the same contradiction. Hence the
sequence is precompact in $X$. The proof is complete.
\end{proof}
The compactness criterion depends on the norm of the weak limit.
No sign condition on $\lambda$ or on $J_\mu(u)$ is needed.

\section{Local minimization and mountain-pass geometry}\label{sec:geometry}
Define
\begin{equation}\label{eq:alphamu}
\mathcal B_\mu=\{u\in M:\norm{\Delta u}_2^2<S_2^{N/4}\mu^{-(N-4)/4}\}.
\end{equation}
If $A=\norm{\Delta u}_2^2\le S_2^{N/4}\mu^{-(N-4)/4}$, it holds $\mu\norm u_q^q\le\mu S_2^{-q/2}A^{q/2}\le A$, which implies that
\begin{equation}\label{eq:smallball}
J_\mu(u)=\frac12\norm{\Delta u}_2^2-\frac\mu q\norm u_q^q\ge\frac2N A\ge\frac2N\Lambda_1.
\end{equation}
On $A=S_2^{N/4}\mu^{-(N-4)/4}$, one has $J_\mu(u)\ge \frac2N S_2^{N/4}\mu^{-(N-4)/4}$.

\begin{proposition}\label{prop:local}
If $\frac2N S_2^{N/4}\mu^{-(N-4)/4}>\Lambda_1/2$, the infimum
\begin{equation}\label{eq:locallevel}
m(\mu)=\inf\{J_\mu(u):u\in M,\ \norm{\Delta u}_2^2\le S_2^{N/4}\mu^{-(N-4)/4}\}
\end{equation}
is attained at $u_\mu^\loc\in\mathcal B_\mu$. Moreover,
\begin{equation}\label{eq:localbounds}
\frac2N\Lambda_1\le m(\mu)<\frac{\Lambda_1}{2}<\frac2N S_2^{N/4}\mu^{-(N-4)/4},
\end{equation}
and $u_\mu^\loc$ is a local minimizer and a critical point of $J_\mu|_M$.
\end{proposition}
\begin{proof}
The infimum defining $\Lambda_1$ is attained at some $\phi_1\in M$ by
compactness into $L^2$. The assumption
gives $S_2^{N/4}\mu^{-(N-4)/4}>N\Lambda_1/4>\Lambda_1$, and
\[
J_\mu(\phi_1)=\frac{\Lambda_1}{2}-\frac\mu q\norm{\phi_1}_q^q
<\frac{\Lambda_1}{2}<\frac2N S_2^{N/4}\mu^{-(N-4)/4}.
\]
Together with \eqref{eq:smallball}, this proves the level bounds.

Let $(v_n)$ minimize on the closed set. We know that it is bounded in $X$.
Pass to a subsequence with $v_n\weak v$ in $X$ and $v_n\to v$ in $L^2$
and almost everywhere. Then $v\in M$ and
$\norm{\Delta v}_2^2\le S_2^{N/4}\mu^{-(N-4)/4}$. With $z_n=v_n-v$, select limits
$A=\lim\norm{\Delta z_n}_2^2$ and $B=\lim\norm{z_n}_q^q$.
Norm splitting gives $0\le A\le S_2^{N/4}\mu^{-(N-4)/4}$, and the Sobolev inequality
gives $\mu B\le\mu S_2^{-q/2}A^{q/2}\le A$. Therefore
\[
m(\mu)=J_\mu(v)+\frac A2-\frac{\mu B}{q}
\ge m(\mu)+\frac2N A.
\]
Thus $A=0$ and $v_n\to v$ in $X$. The boundary energy bound and
$m(\mu)<\frac2N S_2^{N/4}\mu^{-(N-4)/4}$ imply $v\in\mathcal B_\mu$. A relative neighborhood of
$v$ in $M$ lies in the cap; hence $v$ is a local minimizer.
The Lagrange multiplier rule yields \eqref{eq:normalized}.
\end{proof}

\subsection{Cut-off bubbles and endpoints}
Choose $x_0\in\Omega$ and $0\le\eta\in C_c^\infty(\Omega)$ with
$\eta=1$ on $B_R(x_0)\Subset\Omega$. Let
\begin{equation}\label{eq:bubble}
\begin{gathered}
W_\eps(x)=\kappa_N\left(\frac\eps{\eps^2+|x-x_0|^2}\right)^{(N-4)/2},
\qquad w_\eps=\eta W_\eps,\qquad v_\eps=\frac{w_\eps}{\norm{w_\eps}_2},\\
\kappa_N=[N(N-4)(N^2-4)]^{(N-4)/8}.
\end{gathered}
\end{equation}
These are the usual Sobolev optimizers before truncation. In particular,
\begin{equation}\label{eq:bubblenorms}
\norm{\Delta W_\eps}_{L^2(\R^N)}^2
=\norm{W_\eps}_{L^q(\R^N)}^q=S_2^{N/4}.
\end{equation}
The functions $W_\eps$ are the Sobolev optimizers \cite{Sw}, normalized
by $\Delta^2W_\eps=W_\eps^{q-1}$. The norm identities follow by
testing this equation and using the sharp Sobolev quotient.

\begin{lemma}\label{lem:bubble}
For a fixed admissible cutoff, set
$$D_\eta=\int_{\R^N}|\Delta((1-\eta)f)|^2\dd,\ \ \ L_\eta=\int_\Omega\eta^2f^2\dd,$$ with
$f=|x-x_0|^{4-N}$. For every $N\ge5$,
\begin{align}
A_\eps:=\norm{\Delta w_\eps}_2^2
 &=S_2^{N/4}+\kappa_N^2D_\eta\eps^{N-4}+o(\eps^{N-4}),\label{eq:bubbleA}\\
B_\eps:=\norm{w_\eps}_q^q&=S_2^{N/4}+O(\eps^N),\label{eq:bubbleB}\\
L_\eps:=\norm{w_\eps}_2^2&=
\begin{cases}
b_N\eps^4(1+o(1)),&N>8,\\
b_8\eps^4|\log\eps|+O(\eps^4),&N=8,\\
\kappa_N^2L_\eta\eps^{N-4}+o(\eps^{N-4}),&5\le N\le7,
\end{cases}\label{eq:bubbleL}
\end{align}
where $b_N,b_8>0$.
The map $\eps\mapsto v_\eps$ is continuous into $X$ for $\eps>0$.
\end{lemma}
\begin{proof}
Without loss of generality, we suppose $x_0=0$ and choose $0<r_0<R_0$ such that
$\eta=1$ on $B_{r_0}$ and $\operatorname{supp}\eta\subset B_{R_0}$.
Set $H_\eps=(1-\eta)W_\eps$. Since $H_\eps\in D^{2,2}(\R^N)$,
testing the bubble equation gives
\begin{equation}\label{eq:crossbubble}
\int_{\R^N}\Delta W_\eps\Delta H_\eps\dd
=\int_{\R^N}(1-\eta)W_\eps^q\dd=O(\eps^N).
\end{equation}
Indeed, $W_\eps^q\le C\eps^N|x|^{-2N}$ on $|x|\ge r_0$. Moreover,
\begin{equation}\label{eq:Hconvergence}
\frac{H_\eps}{\kappa_N\eps^{(N-4)/2}}
=(1-\eta)(\eps^2+|x|^2)^{-(N-4)/2}
\longrightarrow (1-\eta)|x|^{4-N}
\quad\text{in }D^{2,2}(\R^N).
\end{equation}
Indeed, since both functions vanish on $B_{r_0}$,
we have
\[
\Delta\bigl[(1-\eta)(\eps^2+|x|^2)^{-(N-4)/2}\bigr]
\longrightarrow
\Delta\bigl[(1-\eta)|x|^{4-N}\bigr]
\quad\text{a.e. in }\R^N.
\]
Moreover, for $0<\eps\le 1$ and $|x|\ge r_0$,
\[
\left|
\Delta\bigl[(1-\eta)(\eps^2+|x|^2)^{-(N-4)/2}\bigr]
\right|
\le C(1+|x|)^{2-N}.
\]
Since $N>4$,
$
\int_{\{|x|\ge r_0\}}
(1+|x|)^{4-2N}\,dx<\infty.
$
The same bound holds for the limiting function. Thus \eqref{eq:Hconvergence} follows by dominated convergence and cutoff approximation. Expanding the square and using \eqref{eq:crossbubble}, we obtain
\begin{align*}
A_\eps
&=S_2^{N/4}-2\int_{\R^N}\Delta W_\eps\Delta H_\eps\dd
+\norm{\Delta H_\eps}_2^2\\
&=S_2^{N/4}+\kappa_N^2D_\eta\eps^{N-4}+o(\eps^{N-4}).
\end{align*}
Also,
\[
0\le S_2^{N/4}-B_\eps
=\int_{\R^N}(1-\eta^q)W_\eps^q\dd
\le\int_{|x|\ge r_0}W_\eps^q\dd=O(\eps^N).
\]

For the mass, a change of variables yields
\begin{equation}\label{eq:Lradial}
\begin{aligned}
L_\eps
&=\kappa_N^2\eps^{N-4}\int_\Omega
\eta(x)^2(\eps^2+|x|^2)^{-(N-4)}\dd\\
&=\kappa_N^2\eps^4\int_{\R^N}
\eta(\eps y)^2(1+|y|^2)^{-(N-4)}\,dy.
\end{aligned}
\end{equation}
If $N>8$, the last integrand is dominated by an integrable function,
and hence
\[
L_\eps=b_N\eps^4(1+o(1)),\qquad
b_N=\kappa_N^2\int_{\R^N}(1+|y|^2)^{-(N-4)}dy.
\]
For $N=8$, the integral over $|y|<r_0/\eps$ equals
\[
|\mathbb S^7|\int_0^{r_0/\eps}\frac{r^7}{(1+r^2)^4}\,dr
=|\mathbb S^7|\log(1/\eps)+O(1),
\]
since $r^7(1+r^2)^{-4}=r^{-1}+O(r^{-3})$ for $r\ge1$.
The remaining part is $O(1)$, being bounded by a constant times
$\int_{r_0/\eps}^{R_0/\eps}r^{-1}dr$.
This gives the stated formula with $b_8=\kappa_8^2|\mathbb S^7|$.

If $5\le N\le7$, use instead the first identity in
\eqref{eq:Lradial}. Its normalized integrand converges almost everywhere
to $\eta^2|x|^{8-2N}$ and is bounded by this function. It is integrable
near zero because $\int_0^{r_0}r^{7-N}dr<\infty$. Therefore
$L_\eps/(\kappa_N^2\eps^{N-4})\to L_\eta$.
The positivity of $L_\eta$ follows from $\eta=1$ near zero. Finally, continuity at positive scales follows from the explicit formula on a fixed compact support and $L_\eps>0$.
\end{proof}

Fix $\eps_0\in(0,1)$ sufficiently small and put $v_0=v_{\eps_0}$.
Choose $\mu_*>0$ so small that, for all $0<\mu<\mu_*$,
\begin{equation}\label{eq:start}
\norm{\Delta v_0}_2^2<S_2^{N/4}\mu^{-(N-4)/4},\qquad
J_\mu(v_0)<\frac2N S_2^{N/4}\mu^{-(N-4)/4},\qquad \frac2N S_2^{N/4}\mu^{-(N-4)/4}>\Lambda_1/2.
\end{equation}
This is possible because $S_2^{N/4}\mu^{-(N-4)/4},\frac2N S_2^{N/4}\mu^{-(N-4)/4}\to\infty$.
For a compact interval $I=[\mu_-,\mu_+]\Subset(0,\mu_*)$, choose
$\eps_I\in(0,\eps_0)$ so small that for every $\mu\in I$ and every
$0<\eps\le\eps_I$,
\begin{equation}\label{eq:tail}
\norm{\Delta v_\eps}_2^2>S_2^{N/4}\mu^{-(N-4)/4},\qquad J_\mu(v_\eps)<0.
\end{equation}
Indeed, $A_\eps/L_\eps\to\infty$, while
\[
J_\mu(v_\eps)=\frac{A_\eps}{2L_\eps}
-\frac{\mu B_\eps}{qL_\eps^{q/2}}
\le\frac C{L_\eps}-\frac{c\mu_-}{L_\eps^{q/2}}\to-\infty
\]
uniformly for $\mu\in I$.

\Needspace{12\baselineskip}
\begin{lemma}\label{lem:paths}
Let $I=[\mu_-,\mu_+]\Subset(0,\mu_*)$, and choose $\eps_I$ as in
\eqref{eq:tail}. Set $v_1=v_{\eps_I}$ and
\begin{equation}\label{eq:mountain}
\begin{gathered}
\Gamma_I=\{\gamma\in C([0,1],M):\gamma(0)=v_0,\ \gamma(1)=v_1\},\\
d_I(\mu)=\inf_{\gamma\in\Gamma_I}\max_{t\in[0,1]}J_\mu(\gamma(t)),
\qquad\mu\in I.
\end{gathered}
\end{equation}
Then:
\begin{enumerate}
\item $\Gamma_I\ne\varnothing$ and
\begin{equation}\label{eq:pathlower}
\frac2N S_2^{N/4}\mu^{-(N-4)/4}\le d_I(\mu)<\infty
\qquad(\mu\in I).
\end{equation}
\item Replacing $\eps_I$ by any smaller positive scale satisfying
\eqref{eq:tail} does not change $d_I$. More generally, for two compact
parameter intervals $I,J\Subset(0,\mu_*)$,
\[
d_I(\mu)=d_J(\mu)\qquad(\mu\in I\cap J).
\]
These values define a function $d:(0,\mu_*)\to\R$.
\item The endpoint gap is uniform on $I$:
\begin{equation}\label{eq:uniformpathgap}
\inf_{\mu\in I}
\left[d_I(\mu)-\max\{J_\mu(v_0),J_\mu(v_1)\}\right]>0.
\end{equation}
\end{enumerate}
\end{lemma}
\begin{proof}
\textup{(i)} The path
$\gamma_0(t)=v_{(1-t)\eps_0+t\eps_I}$ belongs to $\Gamma_I$ and $\gamma_0(t)$ for $t\in[0,1]$ is compact in $X$, so its energy maximum is finite.
For any $\gamma\in\Gamma_I$, \eqref{eq:start} and \eqref{eq:tail} give
\[
\norm{\Delta\gamma(0)}_2^2
<S_2^{N/4}\mu^{-(N-4)/4}
<\norm{\Delta\gamma(1)}_2^2.
\]
By continuity there is $t_\gamma\in(0,1)$ such that
$\norm{\Delta\gamma(t_\gamma)}_2^2=S_2^{N/4}\mu^{-(N-4)/4}$.
At this point, Sobolev's inequality yields
\begin{align*}
\mu\norm{\gamma(t_\gamma)}_q^q
&\le\mu S_2^{-q/2}
\left(S_2^{N/4}\mu^{-(N-4)/4}\right)^{q/2}\\
&=S_2^{N/4}\mu^{-(N-4)/4},
\end{align*}
and hence
\[
\max_tJ_\mu(\gamma(t))
\ge J_\mu(\gamma(t_\gamma))
\ge\left(\frac12-\frac1q\right)
S_2^{N/4}\mu^{-(N-4)/4}.
\]
Taking the infimum over $\Gamma_I$ proves \eqref{eq:pathlower}.

\textup{(ii)} Fix $\mu\in I\cap J$. Suppose first that
$\eps_J\le\eps_I$, and define
\[
\zeta(t)=v_{(1-t)\eps_I+t\eps_J},\qquad 0\le t\le1.
\]
Every scale on this path is at most $\eps_I$, so \eqref{eq:tail} gives
\begin{equation}\label{eq:connectornegative}
\max_{t\in[0,1]}J_\mu(\zeta(t))<0.
\end{equation}
For $\gamma\in\Gamma_I$, let
\[
(\gamma*\zeta)(t)=
\begin{cases}
\gamma(2t),&0\le t\le1/2,\\
\zeta(2t-1),&1/2\le t\le1.
\end{cases}
\]
Then $\gamma*\zeta\in\Gamma_J$, and by \textup{(i)} and
\eqref{eq:connectornegative},
\begin{align*}
\max_tJ_\mu((\gamma*\zeta)(t))
&=\max\left\{\max_tJ_\mu(\gamma(t)),
              \max_tJ_\mu(\zeta(t))\right\}\\
&=\max_tJ_\mu(\gamma(t)).
\end{align*}
Therefore $d_J(\mu)\le d_I(\mu)$. Applying the same argument to
paths in $\Gamma_J$ and the reversed connector $\zeta(1-t)$ gives
$d_I(\mu)\le d_J(\mu)$. The case $\eps_I\le\eps_J$ follows by
interchanging the intervals. The assertion for two endpoint scales
on the same interval is the same calculation.

\textup{(iii)} By \textup{(i)},
\begin{align*}
&d_I(\mu)-\max\{J_\mu(v_0),J_\mu(v_1)\}\\
&\qquad\ge\frac2N S_2^{N/4}\mu^{-(N-4)/4}
-\max\{J_\mu(v_0),J_\mu(v_1)\}>0.
\end{align*}
The last expression is continuous in $\mu\in I$ and positive by
\eqref{eq:start} and \eqref{eq:tail}. Its minimum on $I$ is positive,
which proves \eqref{eq:uniformpathgap}.
\end{proof}

For the physical mass $\rho=\mu^{(N-4)/8}$, define
\begin{equation}\label{eq:crho}
c_\rho=\rho^2d(\rho^{q-2}).
\end{equation}
Multiplication by $\rho$ is a homeomorphism $M\to S_\rho$.
Thus \eqref{eq:conversion} identifies this as the mountain-pass value
for the corresponding scaled endpoints on $S_\rho$.

\section{Construction and compactness of mountain-pass solutions}\label{sec:construction}
\subsection{A quantitative monotonicity argument}\label{sec:monotonicity}
We shall use the following form of Ghoussoub's localization theorem
\cite[Theorem 4.5 and Remark 4.10]{G}; see also \cite[Lemma 2.1]{PVY}.

\begin{proposition}[Localized minimax principle]
\label{prop:localized}
Let $\mathcal M$ be a complete Hilbert manifold, $\Phi\in C^1(\mathcal M)$,
and $\Gamma$ the class of paths joining two fixed endpoints.
Suppose the minimax value $c$ exceeds both endpoint values by a fixed
positive gap. There is $\bar\delta>0$, depending on this gap, with the
following property. For $0<\epsilon\le\delta<\bar\delta$, put
$\Gamma_\epsilon=\{\gamma\in\Gamma:\max\Phi(\gamma)\le c+\epsilon\}$.
If a closed set $F$ meets every path in $\Gamma_\epsilon$ and
$\inf_F\Phi\ge c-\delta$, then for any
$\gamma_\epsilon\in\Gamma_\epsilon$ there is $v\in\mathcal M$ such that
\begin{equation}\label{eq:localized}
|\Phi(v)-c|+\norm{d\Phi(v)}+\dist(v,F)
+\dist(v,\gamma_\epsilon([0,1]))\le9\delta+26\sqrt\delta.
\end{equation}
Distances are manifold distances, which dominate ambient Hilbert
distances in our application.
\end{proposition}

\begin{lemma}\label{lem:mono}
The function $d$ is nonincreasing and left-continuous. At every
differentiability point $\mu\in(0,\mu_*)$, there is $(u_n)\subset M$
with
\begin{equation}\label{eq:PSmono}
J_\mu(u_n)\to d(\mu),\qquad
\norm{(J_\mu|_M)'(u_n)}\to0,
\end{equation}
and
\begin{equation}\label{eq:quantitative}
\norm{\Delta u_n}_2^2
\le2d(\mu)-2\mu d'(\mu)+6\mu+o_n(1).
\end{equation}
The error $o_n(1)$ is for fixed $\mu$; uniformity in $\mu$ is not asserted.
\end{lemma}
\begin{proof}
Since $J_\nu(u)\ge J_\mu(u)$ for $\nu<\mu$, $d$ is nonincreasing.
Endpoint compatibility makes this global. If $\mu_n\uparrow\mu$, choose
a path $\gamma$ with $\max J_\mu(\gamma)\le d(\mu)+\epsilon$.
$\gamma(t)$ for $t\in[0,1]$ is compact in $X$, and hence $J_{\mu_n}\to J_\mu$ uniformly
there. Therefore
\[
d(\mu)\le\liminf_n d(\mu_n)\le\limsup_n d(\mu_n)\le d(\mu)+\epsilon.
\]
Let $\epsilon\downarrow0$. Monotone functions are differentiable almost
everywhere, with finite derivative there.

Fix a differentiability point $\mu$ and choose $\mu_n\uparrow\mu$
in the same parameter neighborhood. Put $e_n=\mu-\mu_n$ and
$\delta_n=d(\mu_n)-d(\mu)+e_n$. Define
\begin{align*}
\Gamma_n&=\{\gamma\in\Gamma_I:\max_tJ_{\mu_n}(\gamma(t))
\le d(\mu_n)+e_n\},\\
F_n&=\cl_X\{z=\gamma(t):\gamma\in\Gamma_n,\
J_\mu(z)\ge d(\mu)-e_n\}.
\end{align*}
Every path $\gamma\in\Gamma_n$ there holds $\gamma(t)\cap F_n\neq \emptyset$ since its $J_\mu$ maximum is at
least $d(\mu)$. The set $F_n\subset M$ is closed, and by continuity
every $z\in F_n$ satisfies
\begin{equation}\label{eq:Fbounds}
J_{\mu_n}(z)\le d(\mu_n)+e_n,\qquad J_\mu(z)\ge d(\mu)-e_n.
\end{equation}
Since $J_{\mu_n}(z)\ge J_\mu(z)\geq d(\mu)-e_n= d(\mu_n)-\delta_n$ for $z\in F_n$, we have
$\inf_{F_n}J_{\mu_n}\geq  d(\mu_n)-\delta_n.$
Moreover $0<e_n\le\delta_n\to0$. The uniform endpoint gap in
Lemma~\ref{lem:paths}\textup{(iii)} permits Proposition~\ref{prop:localized} to be
applied to $J_{\mu_n}|_M$ with the same admissible gap for large $n$.
It gives
\begin{equation}\label{eq:approxPS}
\begin{gathered}
J_{\mu_n}(u_n)-d(\mu_n)\to0,\qquad
\norm{(J_{\mu_n}|_M)'(u_n)}\to0,\\
\dist_X(u_n,F_n)\to0.
\end{gathered}
\end{equation}
For every $z\in F_n$, \eqref{eq:Fbounds} implies
\begin{equation}\label{eq:quotient}
\frac1q\norm z_q^q
=\frac{J_{\mu_n}(z)-J_\mu(z)}{\mu-\mu_n}
\le\frac{d(\mu_n)-d(\mu)}{\mu-\mu_n}+2
\le-d'(\mu)+3
\end{equation}
for all sufficiently large $n$. Consequently
\begin{align*}
\frac12\norm{\Delta z}_2^2
&=J_\mu(z)+\frac\mu q\norm z_q^q\\
&\le d(\mu_n)+e_n-\mu d'(\mu)+3\mu.
\end{align*}
This bound is uniform over $F_n$. Choose $z_n\in F_n$ with
$\norm{u_n-z_n}_X\to0$; boundedness of $z_n$ and
$d(\mu_n)\to d(\mu)$ give \eqref{eq:quantitative}. Thus $u_n$ is bounded.

On this bounded sequence, the Sobolev embedding gives
\[
|J_\mu(u_n)-J_{\mu_n}(u_n)|\le Ce_n,\qquad
\norm{J_\mu'(u_n)-J_{\mu_n}'(u_n)}_{X^*}\le Ce_n.
\]
Restriction to tangent spaces decreases the dual norm. Combining these
estimates with \eqref{eq:approxPS} proves \eqref{eq:PSmono}.
\end{proof}

\begin{remark}
The estimate \eqref{eq:quantitative} is the analogue of
\cite[Lemma 2.2(4)]{PVY}. Its derivative term will be used in
Lemma~\ref{lem:strong}. All sequence limits in \eqref{eq:PSmono}--\eqref{eq:quantitative}
are taken at a fixed parameter $\mu$.
\end{remark}

\subsection{An upper bound for the mountain-pass level}\label{sec:upper}
\begin{lemma}\label{lem:upper}
There exists $\mu_*>0$ such that, for $0<\mu<\mu_*$,
\begin{equation}\label{eq:g}
\frac2N S_2^{N/4}\mu^{-(N-4)/4}\le d(\mu)\le
 g(\mu):=\frac2N S_2^{N/4}\mu^{-(N-4)/4}+h(\mu),
\end{equation}
where
\begin{equation}\label{eq:h}
h(\mu)=
\begin{cases}
C\mu^{(N-4)(N-8)/16},&N>8,\\
C|\log\mu|^{-1},&N=8,\\
C_0,&5\le N\le7.
\end{cases}
\end{equation}
Here $C>0$ is independent of $\mu$. For $5\le N\le7$, any
$C_0>D_\eta/(2L_\eta)$ is admissible. More precisely,
\begin{equation}\label{eq:upperasymptotic}
d(\mu)\le\frac2N S_2^{N/4}\mu^{-(N-4)/4}
+\frac{D_\eta}{2L_\eta}+o(1)
\qquad\text{as }\mu\downarrow0.
\end{equation}
\end{lemma}
\begin{proof}
The lower bound follows from Lemma~\ref{lem:paths}. For the upper bound,
use the path $t\mapsto v_{(1-t)\eps_0+t\eps_I}$, with its endpoint chosen
as in \eqref{eq:tail}. Thus
\begin{equation}\label{eq:pathupper}
d(\mu)\le\sup_{0<\eps\le\eps_0}J_\mu(v_\eps).
\end{equation}
Set $\rho=\mu^{(N-4)/8}$ and $s=\rho/L_\eps^{1/2}$. Since
$\rho^2\mu=\rho^q$, we have
\begin{equation}\label{eq:physicalpath}
E(\rho v_\eps)=\rho^2J_\mu(v_\eps)
=\frac12A_\eps s^2-\frac1qB_\eps s^q.
\end{equation}
We may suppose $S_2^{N/4}/2\le A_\eps,B_\eps\le2S_2^{N/4}$ on
$(0,\eps_0]$. If $E(\rho v_\eps)\ge S_2^{N/4}/N,$ then
\[S_2^{N/4}s^2\ge\frac12A_\eps s^2\ge \frac12A_\eps s^2-\frac1qB_\eps s^q
\ge S_2^{N/4}/N,
\]
\[S_2^{N/4}\left(s^2-\frac{s^q}{2q}\right)
\ge\frac12A_\eps s^2-\frac1qB_\eps s^q
\ge S_2^{N/4}/N>0,
\]
and hence
\begin{equation}\label{eq:sbound}
N^{-1/2}\le s\le(2q)^{1/(q-2)},\qquad L_\eps\asymp\rho^2.
\end{equation}
All scales satisfying this energy inequality tend to zero uniformly as
$\rho\to0$. Otherwise, some such scales would remain in a compact
subinterval of $(0,\eps_0]$, on which $L_\eps$ has a positive minimum,
contradicting \eqref{eq:sbound}.

The function $s\mapsto s^2/2-s^q/q$ attains its maximum at $s=1$.
It follows from Lemma~\ref{lem:bubble} and \eqref{eq:sbound} that
\begin{equation}\label{eq:coarseE}
E(\rho v_\eps)
\le\frac2N S_2^{N/4}+C\eps^{N-4}.
\end{equation}
We distinguish the three dimensional ranges.

If $N>8$, \eqref{eq:bubbleL} and \eqref{eq:sbound} give
$\eps^4\asymp\rho^2=\mu^{(N-4)/4}$, so that
\[
\rho^{-2}\eps^{N-4}
\le C\mu^{-(N-4)/4+(N-4)^2/16}
=C\mu^{(N-4)(N-8)/16}.
\]
Dividing \eqref{eq:coarseE} by $\rho^2$ proves the desired bound.

Let $N=8$. Then $\eps^4|\log\eps|\asymp\mu$.
Taking logarithms yields
\[
|\log\mu|=4|\log\eps|-\log|\log\eps|+O(1),
\]
so $|\log\eps|\asymp|\log\mu|$, uniformly over these scales.
Consequently $\rho^{-2}\eps^4\le C|\log\mu|^{-1}$, and the result
again follows from \eqref{eq:coarseE}.

For $5\le N\le7$, the precise leading coefficient is needed. We have
$\eps^{N-4}\asymp\rho^2$ and $\eps\asymp\mu^{1/4}$. By the
expansions of $A_\eps$ and $L_\eps$,
\begin{equation}\label{eq:exactcorrection}
\frac12(A_\eps-S_2^{N/4})s^2
=\rho^2\frac{A_\eps-S_2^{N/4}}{2L_\eps}
=\rho^2\left(\frac{D_\eta}{2L_\eta}+o(1)\right).
\end{equation}
The error is uniform, since all the relevant scales tend to zero and
$L_\eps/\eps^{N-4}\to\kappa_N^2L_\eta>0$.
The critical-norm correction satisfies
\[
\frac{|B_\eps-S_2^{N/4}|s^q}{\rho^2}
\le C\frac{\eps^N}{\rho^2}
\le C\eps^4=o(\mu).
\]
Using \eqref{eq:physicalpath}, the scalar maximum and
\eqref{eq:exactcorrection}, we conclude that
\[
E(\rho v_\eps)\le\frac2N S_2^{N/4}
+\rho^2\left(\frac{D_\eta}{2L_\eta}+o(1)\right).
\]
After division by $\rho^2$, this proves \eqref{eq:upperasymptotic}
on the scales with $E(\rho v_\eps)\ge S_2^{N/4}/N$.

For all remaining scales,
$J_\mu(v_\eps)<S_2^{N/4}\mu^{-(N-4)/4}/N$, which is below the claimed
upper bound. Thus the estimates hold on the whole path. Taking the
supremum in \eqref{eq:pathupper}, and then decreasing $\mu_*$ if necessary,
proves \eqref{eq:g} in every dimension.
\end{proof}

\subsection{Parameter selection and compactness}\label{sec:selection}
We assume (H) if $5\le N\le7$, and fix $C_0$ as in
\eqref{eq:C0choice}. In that case choose $d_0$ with
\begin{equation}\label{eq:d0}
\frac{NC_0}{2S_2^{N/4}}<d_0<
\frac{N(\Lambda_1-2C_0)}{(N-4)S_2^{N/4}}.
\end{equation}
The interval is nonempty, since $NC_0<2\Lambda_1$.

\begin{lemma}\label{lem:selection}
After decreasing $\mu_*$, there is a measurable set $D\subset(0,\mu_*)$
such that $|D\cap(0,t)|>0$ for every $0<t<\mu_*$, and $d$ is
differentiable at every $\mu\in D$, with
\begin{equation}\label{eq:derivativecomparison}
-d'(\mu)\le(1+\delta(\mu))[-g'(\mu)],
\end{equation}
where
\begin{equation}\label{eq:delta}
\delta(\mu)=
\begin{cases}
\mu^{N(N-4)/32},&N>8,\\
\mu|\log\mu|^{-1/2},&N=8,\\
d_0\mu^{(N-4)/4},&5\le N\le7.
\end{cases}
\end{equation}
\end{lemma}
\begin{proof}
The explicit upper bound satisfies
\[
g'(\mu)=-\frac{N-4}{2N}S_2^{N/4}\mu^{-N/4}+h'(\mu).
\]
More precisely,
\[
\frac{h'(\mu)}{\frac{N-4}{2N}S_2^{N/4}\mu^{-N/4}}
=
\begin{cases}
\dfrac{CN(N-8)}{8S_2^{N/4}}\mu^{(N-4)^2/16},&N>8,\\[4pt]
\dfrac{4C}{S_2^2}\mu|\log\mu|^{-2},&N=8,\\[4pt]
0,&N=5,6,7,
\end{cases}
\longrightarrow0.
\]
Thus, after decreasing $\mu_*$, $g'<0$ on $(0,\mu_*)$.
The function $\delta$ in \eqref{eq:delta} is increasing in each case.

Let $D$ be the set of finite differentiability points of $d$ at which
\eqref{eq:derivativecomparison} holds. It is measurable. Suppose that
$|D\cap(0,t_0)|=0$ for some $t_0\in(0,\mu_*)$. Since $d$ is
nonincreasing, it is differentiable almost everywhere, and
\begin{equation}\label{eq:badderivative}
-d'(s)>(1+\delta(s))[-g'(s)]
\quad\text{for almost every }s\in(0,t_0).
\end{equation}
Set $\gamma=2$ for $N\ge8$. For $N=5,6,7$, \eqref{eq:d0} gives
\[
0<\frac{NC_0}{2d_0S_2^{N/4}}<1.
\]
Choose
\[
\gamma>
\left(1-\frac{NC_0}{2d_0S_2^{N/4}}\right)^{-4/(N-4)}.
\]
Then $\gamma>1$ and
\begin{equation}\label{eq:gamma}
\frac{2d_0}{N}S_2^{N/4}
(1-\gamma^{-(N-4)/4})>C_0.
\end{equation}

For $0<t<\gamma t<t_0$, the monotonicity of $d$ implies
\begin{equation}\label{eq:monotonevariation}
d(t)-d(\gamma t)\ge\int_t^{\gamma t}-d'(s)\,ds.
\end{equation}
Indeed, if $0<e<(\gamma-1)t$, then
\begin{align*}
0\le\frac1e\int_t^{\gamma t-e}[d(s)-d(s+e)]\,ds
&=\frac1e\int_t^{t+e}d(s)\,ds
-\frac1e\int_{\gamma t-e}^{\gamma t}d(s)\,ds\\
&\le d(t)-d(\gamma t).
\end{align*}
Extend the nonnegative difference quotient by zero on
$(\gamma t-e,\gamma t)$ and apply Fatou's lemma as $e\downarrow0$.
This proves \eqref{eq:monotonevariation}.
Combining \eqref{eq:badderivative} with $-g'>0$ and the monotonicity
of $\delta$, we obtain
\begin{align}\label{eq:variationselection}
d(t)-d(\gamma t)
&\ge\int_t^{\gamma t}(1+\delta(s))[-g'(s)]\,ds\notag\\
&\ge(1+\delta(t))[g(t)-g(\gamma t)].
\end{align}
On the other hand, \eqref{eq:g} yields
\[
d(t)-d(\gamma t)
\le g(t)-[g(\gamma t)-h(\gamma t)].
\]
Subtracting the common difference from these two estimates gives
\begin{equation}\label{eq:selectioncontradiction}
\delta(t)[g(t)-g(\gamma t)]\le h(\gamma t).
\end{equation}
We show that this is impossible in each dimensional range.

If $N>8$, direct subtraction gives
\begin{align*}
g(t)-g(2t)
={}&\frac2N S_2^{N/4}(1-2^{-(N-4)/4})t^{-(N-4)/4}\\
&-C\big(2^{(N-4)(N-8)/16}-1\big)t^{(N-4)(N-8)/16}\\
\ge{}&\frac1N S_2^{N/4}(1-2^{-(N-4)/4})t^{-(N-4)/4}
\end{align*}
for all sufficiently small $t$. In fact, the ratio of the magnitude of
the second term to the first is $O(t^{(N-4)^2/16})\to0$.
Using \eqref{eq:delta} and \eqref{eq:selectioncontradiction}, we infer
\begin{align*}
&\frac1N S_2^{N/4}(1-2^{-(N-4)/4})
 t^{(N-4)(N-8)/32}\\
&\qquad\le\delta(t)[g(t)-g(2t)]
\le C2^{(N-4)(N-8)/16}t^{(N-4)(N-8)/16}.
\end{align*}
which gives a contradiction for small $t$.

If $N=8$, then for sufficiently small $t$,
\begin{align*}
g(t)-g(2t)
&=\frac{S_2^2}{8t}
-\frac{C\log2}{|\log t|(|\log t|-\log2)}\\
&\ge\frac{S_2^2}{16t},
\qquad
h(2t)=\frac{C}{|\log t|-\log2}\le\frac{2C}{|\log t|}.
\end{align*}
Thus \eqref{eq:selectioncontradiction} gives
\[
\frac{S_2^2}{16}|\log t|^{-1/2}
\le\delta(t)[g(t)-g(2t)]
\le\frac{2C}{|\log t|},
\]
and consequently again a contradiction.

Finally, for $N=5,6,7$, $h=C_0$ and
\[
g(t)-g(\gamma t)
=\frac2N S_2^{N/4}t^{-(N-4)/4}
(1-\gamma^{-(N-4)/4}).
\]
Equations \eqref{eq:gamma} and \eqref{eq:selectioncontradiction} then
imply the impossible chain
\[
C_0<\frac{2d_0}{N}S_2^{N/4}(1-\gamma^{-(N-4)/4})
=\delta(t)[g(t)-g(\gamma t)]
\le h(\gamma t)=C_0.
\]
Therefore $|D\cap(0,t)|>0$ for every $0<t<\mu_*$.
\end{proof}

\begin{remark}
The preceding argument gives positive measure in every interval adjacent
to zero. It does not provide an almost-everywhere derivative bound on
the whole parameter interval.
\end{remark}

\begin{lemma}\label{lem:strong}
For a sufficiently small $\mu_*$, the sequence in Lemma~\ref{lem:mono}
has a strongly convergent subsequence at every $\mu\in D$. Its limit
$u_\mu\in M$ satisfies
\begin{equation}\label{eq:selectedbounds}
J_\mu(u_\mu)=d(\mu),\qquad (J_\mu|_M)'(u_\mu)=0,\qquad
\norm{\Delta u_\mu}_2^2
\le S_2^{N/4}\mu^{-(N-4)/4}+\Lambda_1-\sigma,
\end{equation}
where $\sigma>0$ is independent of $\mu$. One may take
$\sigma=\Lambda_1/2$ for $N\ge8$, and
\begin{equation}\label{eq:sigmachoice}
\sigma=\frac12\left(\Lambda_1-2C_0
-\frac{N-4}{N}S_2^{N/4}d_0\right)>0
\quad\text{for }5\le N\le7.
\end{equation}
\end{lemma}
\begin{proof}
Fix $\mu\in D$ and let $(u_n)$ be the sequence in Lemma~\ref{lem:mono}.
By \eqref{eq:g} and \eqref{eq:derivativecomparison},  we have
\begin{align}\label{eq:Rbound}
\limsup_n\norm{\Delta u_n}_2^2
&\le2g(\mu)-2\mu(1+\delta(\mu))g'(\mu)+6\mu\notag\\
&=S_2^{N/4}\mu^{-(N-4)/4}+2h(\mu)
+\frac{N-4}{N}S_2^{N/4}\delta(\mu)\mu^{-(N-4)/4}\notag\\
&\qquad-2\mu(1+\delta(\mu))h'(\mu)+6\mu.
\end{align}
Since $h'\ge0$,
the last derivative term is nonpositive. If $N>8$, the upper remainder in \eqref{eq:Rbound} is bounded by
\[
2C\mu^{(N-4)(N-8)/16}
+\frac{N-4}{N}S_2^{N/4}\mu^{(N-4)(N-8)/32}+6\mu=o(1).
\]
For $N=8$, it is bounded by
\[
\frac{2C}{|\log\mu|}
+\frac{S_2^2}{2|\log\mu|^{1/2}}+6\mu=o(1).
\]
We can therefore decrease $\mu_*$ so that either remainder is at most
$\Lambda_1/2$. If $5\le N\le7$, the remainder is
\[
2C_0+\frac{N-4}{N}S_2^{N/4}d_0+6\mu.
\]
The choice \eqref{eq:d0} makes its constant part strictly smaller than
$\Lambda_1$. With $\sigma$ as in \eqref{eq:sigmachoice}, choose
$\mu_*$ so that $6\mu_*\le\sigma$. Thus in every case
\begin{equation}\label{eq:strictmargin}
\limsup_n\norm{\Delta u_n}_2^2
\le S_2^{N/4}\mu^{-(N-4)/4}+\Lambda_1-\sigma
<S_2^{N/4}\mu^{-(N-4)/4}+\Lambda_1.
\end{equation}

Lemma~\ref{lem:compact} now gives strong convergence along a subsequence.
Indeed, the weak limit has unit $L^2$ norm and hence squared $X$ norm
at least $\Lambda_1$, while any nonzero defect would have squared norm
at least $S_2^{N/4}\mu^{-(N-4)/4}$, contradicting
\eqref{eq:strictmargin}. Strong convergence in $X$ also gives convergence
in $L^q$, so the level and norm bounds pass to the limit. The limiting
weak equation gives $(J_\mu|_M)'(u_\mu)=0$.
\end{proof}
\section{Proofs of the main results}\label{sec:results}
\subsection{Proof of the main theorem}\label{sec:mainproof}
\begin{proposition}\label{prop:closure}
For every $\hat\mu\in\overline D\cap(0,\mu_*)$, there is a constrained
critical point $\hat u\in M$ of $J_{\hat\mu}$ satisfying
\begin{equation}\label{eq:closurebounds}
\begin{gathered}
\frac2N S_2^{N/4}\hat\mu^{-(N-4)/4}
\le J_{\hat\mu}(\hat u)\le g(\hat\mu),\\
\norm{\Delta\hat u}_2^2
\le S_2^{N/4}\hat\mu^{-(N-4)/4}+\Lambda_1-\sigma.
\end{gathered}
\end{equation}
If $\hat\mu\in D$, or if $\mu_j\uparrow\hat\mu$ for some sequence
$(\mu_j)\subset D$, then $\hat u$ can be chosen with
$J_{\hat\mu}(\hat u)=d(\hat\mu)$.
\end{proposition}
\begin{proof}
The case $\hat\mu\in D$ follows from Lemma~\ref{lem:strong}.
If $\hat\mu\in \overline{D}\setminus D$, choose $\mu_j\in D$ with $\mu_j\to\hat\mu$, and let
$u_j=u_{\mu_j}$ be the corresponding critical points. By
\eqref{eq:selectedbounds}, $(u_j)$ is bounded in $X$, since
$\hat\mu>0$. Passing to a subsequence, we have $d(\mu_j)\to\ell$,
where
\[
\frac2N S_2^{N/4}\hat\mu^{-(N-4)/4}\le\ell\le g(\hat\mu).
\]
On this bounded sequence, H\"older's inequality and the Sobolev embedding
give
\begin{align*}
|J_{\hat\mu}(u_j)-J_{\mu_j}(u_j)|
&\le C|\hat\mu-\mu_j|,\\
\norm{J_{\hat\mu}'(u_j)-J_{\mu_j}'(u_j)}_{X^*}
&\le C|\hat\mu-\mu_j|.
\end{align*}
Thus $(u_j)$ is a constrained Palais--Smale sequence for $J_{\hat\mu}$
at level $\ell$. Moreover,
\[
\limsup_j\norm{\Delta u_j}_2^2
\le S_2^{N/4}\hat\mu^{-(N-4)/4}+\Lambda_1-\sigma
<S_2^{N/4}\hat\mu^{-(N-4)/4}+\Lambda_1.
\]
Lemma~\ref{lem:compact} yields strong convergence along a subsequence
to a critical point $\hat u\in M$. Passing to the limit gives
\eqref{eq:closurebounds}. If $\mu_j\uparrow\hat\mu$, the
left-continuity of $d$ implies $\ell=d(\hat\mu)$.
\end{proof}

\begin{proof}[Proof of Theorem~\ref{thm:main}]
\textup{(i)} Choose $\mu_*>0$ small enough that
Proposition~\ref{prop:local} applies for $0<\mu<\mu_*$, and put
$\rho_*=\mu_*^{(N-4)/8}$. For $0<\rho<\rho_*$, set
\[
\mu=\rho^{8/(N-4)},\qquad U_\rho^{\loc}=\rho u_\mu^{\loc}.
\]
Since multiplication by $\rho$ maps $M$ homeomorphically onto $S_\rho$
and $E(\rho u)=\rho^2J_\mu(u)$, the function $U_\rho^{\loc}$ is a
local minimizer of $E|_{S_\rho}$. Scaling the weak equation gives
equation \eqref{eq:original}, with the same multiplier, because
$\mu\rho^{2-q}=1$. Moreover, \eqref{eq:localbounds} implies
\[
\frac2N\Lambda_1\rho^2\le E(U_\rho^{\loc})
<\frac{\Lambda_1\rho^2}{2}.
\]
Decrease $\rho_*$ so that
$\Lambda_1\rho_*^2/2<2S_2^{N/4}/N$. This proves \eqref{eq:localmain}
without hypothesis (H).

\textup{(ii)} Let $D$, $\mu_*$ and $\sigma$ be as in
Lemmas~\ref{lem:selection} and \ref{lem:strong}, with $\mu_*$ decreased
if necessary, and set
\[
\rho_*=\mu_*^{(N-4)/8},\qquad
Q=\{\mu^{(N-4)/8}:\mu\in D\}.
\]
The parameter map is an increasing homeomorphism and is locally
bi-Lipschitz away from zero, so $Q$ is measurable. The change of variables
gives, for $0<t<\rho_*$,
\begin{equation}\label{eq:masssetmeasure}
|Q\cap(0,t)|
=\int_{D\cap(0,t^{8/(N-4)})}
\frac{N-4}{8}\mu^{(N-12)/8}\,d\mu>0.
\end{equation}
If $\rho\in\overline Q\cap(0,\rho_*)$, then
$\mu=\rho^{8/(N-4)}\in\overline D\cap(0,\mu_*)$.
Let $\hat u$ be given by Proposition~\ref{prop:closure} and put
$U_\rho^2=\rho\hat u$. The scaling used in \textup{(i)} gives a solution
of equation \eqref{eq:original}. Multiplying \eqref{eq:closurebounds} by
$\rho^2=\mu^{(N-4)/4}$ yields
\begin{align*}
\frac2N S_2^{N/4}\le E(U_\rho^2)
&\le\frac2N S_2^{N/4}+\rho^2h(\mu),\\
\norm{\Delta U_\rho^2}_2^2
&\le S_2^{N/4}+(\Lambda_1-\sigma)\rho^2.
\end{align*}
In view of \eqref{eq:h},
\[
\rho^2h(\rho^{8/(N-4)})=
\begin{cases}
C\rho^{(N-4)/2},&N>8,\\
\dfrac C2\rho^2/|\log\rho|,&N=8,\\
C_0\rho^2,&N=5,6,7.
\end{cases}
\]
This proves \eqref{eq:secondmain}--\eqref{eq:secondnorm}, after
renaming the logarithmic constant. The inequality
\[
E(U_\rho^{\loc})<\frac2N S_2^{N/4}\le E(U_\rho^2)
\]
implies $U_\rho^2\ne\pm U_\rho^{\loc}$.
Finally, at selected masses or increasing limits of selected masses,
Proposition~\ref{prop:closure} gives $J_\mu(\hat u)=d(\mu)$; hence
$E(U_\rho^2)=\rho^2d(\mu)=c_\rho$ by \eqref{eq:crho}.
\end{proof}

\subsection{Verification of the low-dimensional hypothesis}\label{sec:explicit}
Throughout this section $N=5,6,7$. We construct cutoffs satisfying (H)
and prove Corollary~\ref{cor:geometry}.

\subsubsection{Test profiles and smooth cutoffs}
Let $r=|x|$, $f=r^{4-N}$, and define $G=0$ for $r\ge1$, while
\begin{equation}\label{eq:radialG}
G(r)=
\begin{cases}
r^{-1}-\dfrac32+\dfrac12r^2+(1-r^2)^2,&N=5,\\[2pt]
r^{-2}-2+r^2,&N=6,\\[2pt]
r^{-3}-\dfrac52+\dfrac32r^2,&N=7,
\end{cases}
\qquad 0<r<1.
\end{equation}
Set $H=f-G$. Thus, on $r<1$,
\begin{equation}\label{eq:Hprofile}
H(r)=
\begin{cases}
\dfrac12+\dfrac32r^2-r^4,&N=5,\\[2pt]
2-r^2,&N=6,\\[2pt]
\dfrac52-\dfrac32r^2,&N=7,
\end{cases}
\end{equation}
and $H=f$ for $r\ge1$. Since \(G(1)=G'(1)=0\), the corresponding first-order traces of the piecewise function \(H\) match at \(r=1\). Moreover,
\[
\Delta H=2(4-N)r^{2-N}\qquad (r>1),
\]
and hence, for \(N>4\),
\[
\int_{|x|>1}|\Delta H|^2\,dx
=4(N-4)^2|\mathbb S^{N-1}|
\int_1^\infty r^{3-N}\,dr<\infty.
\]
Together with the local regularity of \(H\), this yields
\(
H\in D^{2,2}(\mathbb R^N).
\)
Indeed, cutoff approximation at infinity gives a squared error of order \(O(R^{4-N})\), while the matching first-order traces exclude any distributional surface term.

To establish \(0\le G\le f\), observe that
\[
\phi(r):=r^{4-N}-\frac{N-2}{2}
+\frac{N-4}{2}r^2
\]
satisfies
\[
\phi(1)=0,\qquad
\phi'(r)=(N-4)(r-r^{3-N})<0
\quad (0<r<1).
\]
Thus \(\phi(r)>0\) for \(0<r<1\). In addition, for \(N=5\),
\[
H(r)=\frac12+\frac32r^2-r^4\ge\frac12,
\]
whereas for \(N=6,7\), the inner tail is bounded below by \(1\). Consequently,
\[
0\le G(r)\le f(r),\qquad 0<r<1.
\]
Furthermore,
\[
\int_{B_1}|G|^2\,dx
\le C\int_0^1r^{7-N}\,dr<\infty,
\qquad N=5,6,7.
\]
Although \(G\notin H_0^2(B_1)\), it can be used to define the limiting cutoff quotient through admissible smooth approximations. Write
\[
D_G=\int_{\R^N}|\Delta H|^2\dd,\qquad
L_G=\int_{B_1}G^2\dd,\qquad \omega=|\mathbb S^{N-1}|.
\]
For $N=6,7$, the inner Laplacian is $-N(N-4)$ and the outer one is
$-2(N-4)r^{2-N}$. Thus
\begin{align}\label{eq:Dbase}
\frac{D_G}{\omega}
&=N^2(N-4)^2\int_0^1r^{N-1}dr
+4(N-4)^2\int_1^\infty r^{3-N}dr\notag\\
&=N(N-4)^2+4(N-4)=(N-4)(N-2)^2.
\end{align}
For $N=5$, let $\psi=(1-r^2)^2$ on $B_1$, extended by zero.
The baseline tail $3/2-r^2/2$ has energy $9\omega$, and
\[
\Delta\psi=-4N+4(N+2)r^2,\qquad \int_{B_1}\Delta\psi\dd=0.
\]
Hence its cross term with the constant inner Laplacian of the baseline
tail vanishes. Furthermore,
\[
\frac1\omega\int_{B_1}|\Delta\psi|^2\dd
=16\left(N-2N+\frac{(N+2)^2}{N+4}\right)
=\frac{64}{N+4}.
\]
Consequently $D_G/\omega=9+64/9=145/9$ for $N=5$.

For the baseline mass, expansion of the square gives
\begin{equation}\label{eq:Lbase}
\begin{aligned}
\int_0^1\left(r^{4-N}-\frac{N-2}{2}+\frac{N-4}{2}r^2\right)^2r^{N-1}dr
={}&\frac1{8-N}+\frac{(N-2)^2}{4N}+\frac{(N-4)^2}{4(N+4)}\\
&-\frac{N-2}{4}+\frac{N-4}{6}
-\frac{(N-2)(N-4)}{2(N+2)}.
\end{aligned}
\end{equation}
This yields $L_G/\omega=1/10$ for $N=6$ and $475/924$ for $N=7$.
For $N=5$, the baseline mass is $17/1260$, and
\begin{align*}
\int_0^1(r^{-1}-\tfrac32+\tfrac12r^2)\psi r^4dr
&=\left(\frac14-\frac26+\frac18\right)
-\frac32\left(\frac15-\frac27+\frac19\right)\\
&\quad+\frac12\left(\frac17-\frac29+\frac1{11}\right)
=\frac{37}{3960},\\
\int_0^1\psi^2r^4dr
&=\frac15-\frac47+\frac69-\frac4{11}+\frac1{13}
=\frac{128}{15015}.
\end{align*}
It follows that
\[
\frac{L_G}{\omega}
=\frac{17}{1260}+2\frac{37}{3960}+\frac{128}{15015}
=\frac{3667}{90090}.
\]
The resulting quotients are listed in Table~\ref{tab:constants}.

\begin{table}[ht]
\centering
\caption{Constants associated with the profiles in \eqref{eq:radialG}.}
\label{tab:constants}
\begin{tabular}{@{}crrrr@{}}
\toprule
$N$&$D_G/\omega$&$L_G/\omega$&$D_G/(2L_G)$&$ND_G/(4L_G)$\\
\midrule
5&$145/9$&$3667/90090$&$725725/3667$&$3628625/7334$\\
6&$32$&$1/10$&$160$&$480$\\
7&$75$&$475/924$&$1386/19$&$4851/19$\\
\bottomrule
\end{tabular}
\end{table}

\begin{lemma}\label{lem:smoothing}
There exist $\eta_h\in C_c^\infty(B_1)$, $0\le\eta_h\le1$, equal to
one near zero, such that
\begin{equation}\label{eq:cutofflimit}
D_{\eta_h}\to D_G,\qquad L_{\eta_h}\to L_G
\quad\text{as }h\downarrow0.
\end{equation}
In particular $D_{\eta_h}/(2L_{\eta_h})\to D_G/(2L_G)$.
\end{lemma}
\begin{proof}
For $0<h<1/8$, choose radial cutoffs $\chi_h,\zeta_h$ taking values
in $[0,1]$ such that
\[
\begin{gathered}
\chi_h=0\text{ on }r\le h,\qquad \chi_h=1\text{ on }r\ge2h,\\
\zeta_h=1\text{ on }r\le1-2h,\qquad
\zeta_h=0\text{ on }r\ge1-h,
\end{gathered}
\]
with first and second derivatives bounded by $Ch^{-1}$ and $Ch^{-2}$.
Define, for $0<r<1$,
\[
\eta_h=\zeta_h\left(1-\chi_h\frac Hf\right),
\]
and extend it by one at zero and by zero outside $B_1$.
Since it is identically one on $B_h$ and vanishes for $r\ge1-h$,
$\eta_h\in C_c^\infty(B_1)$. The inequalities $0\le H\le f$ give
$0\le\eta_h\le1$. Also,
\begin{equation}\label{eq:smoothingpieces}
(1-\eta_h)f-H=(\chi_h-1)H+(1-\zeta_h)G.
\end{equation}
The two terms on the right have disjoint supports.

Near the origin, $H$ and its first two derivatives are bounded. The
product rule therefore yields
\[
\norm{\Delta((\chi_h-1)H)}_2^2
\le C(h^N+h^{N-2}+h^{N-4})\le Ch^{N-4}.
\]
On the boundary shell, $G(1)=G'(1)=0$ implies
$|G|\le Ch^2$, $|\nabla G|\le Ch$ and $|D^2G|\le C$.
Thus the terms in $\Delta((1-\zeta_h)G)$ are bounded uniformly, and
integration over the shell gives
\[
\norm{\Delta((1-\zeta_h)G)}_2^2\le Ch.
\]
The zero clamped traces ensure that no boundary distribution is introduced
by extension. By \eqref{eq:smoothingpieces},
\begin{equation}\label{eq:smoothingenergy}
\norm{\Delta((1-\eta_h)f-H)}_2^2\le C(h^{N-4}+h)\to0.
\end{equation}
Similarly, boundedness of $H$ near zero and $G=O(h^2)$ near the boundary
give
\[
\norm{\eta_hf-G}_{L^2(B_1)}^2\le C(h^N+h^5)\to0.
\]
The two norm convergences imply \eqref{eq:cutofflimit}. Since $L_G>0$,
the quotients converge as well.
\end{proof}

\subsubsection{Eigenvalue bounds on the ball}
The clamped biharmonic Green function on a ball is positive; see
\cite{GR}. Accordingly, the first clamped eigenfunction on $B_1$ may be
chosen positive. Let $e_1$ be such an eigenfunction. If
$\phi\in H_0^2(B_1)\cap C^4(\overline{B_1})$ is positive in $B_1$
and $\Delta^2\phi\ge c\phi$, then
\begin{equation}\label{eq:positivecomparison}
\Lambda_1(B_1)\int_{B_1}e_1\phi\dd
=\int_{B_1}\Delta e_1\Delta\phi\dd
=\int_{B_1}e_1\Delta^2\phi\dd
\ge c\int_{B_1}e_1\phi\dd.
\end{equation}
Thus $\Lambda_1(B_1)\ge c$, with strict inequality if the last
pointwise inequality is strict on a set of positive measure.
We apply this comparison to explicit polynomials. For $k\ge2$,
\begin{equation}\label{eq:bilappoly}
\Delta^2r^{2k}
=2k(2k+N-2)(2k-2)(2k+N-4)r^{2k-4}.
\end{equation}

For $N=7$, choose $\phi=(1-r^2)^2$. Then $0<\phi\le1$ and
$\Delta^2\phi=8N(N+2)=504$, so
\[
\Lambda_1(B_1)\ge504>4851/19=T_7.
\]
For $N=6$, set $x=r^2$ and
$\phi=(1-x)^2(1-x/5)$. Expanding gives
$\phi=1-11x/5+7x^2/5-x^3/5$, whence
\[
\Delta^2\phi=2688/5-384x
\ge(2688/5)(1-x)\ge(2688/5)\phi.
\]
Therefore $\Lambda_1(B_1)\ge2688/5>480=T_6$.

For $N=5$, take
\begin{align*}
\phi&=(1-x)^2(1-3x/5+3x^2/20)\\
&=1-13x/5+47x^2/20-9x^3/10+3x^4/20.
\end{align*}
The last factor is at least $11/20$ on $[0,1]$, so $\phi$ is positive
inside the ball and has clamped traces. By \eqref{eq:bilappoly},
\[
\Delta^2\phi=658-\frac{6804}{5}x+\frac{3564}{5}x^2.
\]
Writing $y=1-x\in[0,1]$, we obtain
\begin{align*}
\Delta^2\phi-500\phi
&=10-\frac{324}{5}y+\frac{2189}{5}y^2-150y^3-75y^4\\
&\ge10-\frac{324}{5}y+\frac{1064}{5}y^2>0.
\end{align*}
The last quadratic has positive leading coefficient and discriminant
$-107824/25$. Hence
\[
\Lambda_1(B_1)>500>3628625/7334=T_5.
\]
All three test functions are smooth, positive in the ball, and vanish
together with their first normal derivative on its boundary.

\subsubsection{Proof of the domain condition}
If $B_R(x_0)\subset\Omega$, define
$\eta_R(x)=\eta_h((x-x_0)/R)$. A change of variables gives
\begin{equation}\label{eq:cutoffscaling}
D_{\eta_R}=R^{4-N}D_{\eta_h},\qquad
L_{\eta_R}=R^{8-N}L_{\eta_h},\qquad
\frac{D_{\eta_R}}{2L_{\eta_R}}
=R^{-4}\frac{D_{\eta_h}}{2L_{\eta_h}}.
\end{equation}

\begin{proof}[Proof of Corollary~\ref{cor:geometry}]
\textup{(i)} Assume \eqref{eq:domaincondition}. Since the inequality
is strict, choose $0<R<R_\Omega$ such that
\[
\Lambda_1(\Omega)R^4>T_N.
\]
By the definition of the inradius, there is $x_0\in\Omega$ such that
$B_R(x_0)\subset\Omega$. The constants in
Table~\ref{tab:constants} satisfy $T_N=ND_G/(4L_G)$; hence
\[
R^{-4}\frac{D_G}{2L_G}<\frac{2\Lambda_1(\Omega)}{N}.
\]
By Lemma~\ref{lem:smoothing},
\[
R^{-4}\frac{D_{\eta_h}}{2L_{\eta_h}}
\longrightarrow R^{-4}\frac{D_G}{2L_G}
\qquad\text{as }h\downarrow0.
\]
Fix $h>0$ sufficiently small that
\[
R^{-4}\frac{D_{\eta_h}}{2L_{\eta_h}}
<\frac{2\Lambda_1(\Omega)}{N}.
\]
The function $\eta_R(x)=\eta_h((x-x_0)/R)$ belongs to
$C_c^\infty(\Omega)$, takes values in $[0,1]$, and is equal to one
near $x_0$. Equation~\eqref{eq:cutoffscaling} now gives
\[
\frac{D_{\eta_R}}{2L_{\eta_R}}
=R^{-4}\frac{D_{\eta_h}}{2L_{\eta_h}}
<\frac{2\Lambda_1(\Omega)}{N},
\]
which is precisely \eqref{eq:cutoffhyp}.

\medskip
\noindent\textup{(ii)} Let $\Omega=B_R(x_0)$. Since $R_\Omega=R$,
the change of variables $x=x_0+Ry$ in the Rayleigh quotient yields
\[
\Lambda_1(B_R(x_0))R_\Omega^4=\Lambda_1(B_1).
\]
The eigenvalue estimates proved above give
\[
\begin{array}{ll}
N=5:&\displaystyle\Lambda_1(B_1)>500>
\frac{3628625}{7334}=T_5,\\[4pt]
N=6:&\displaystyle\Lambda_1(B_1)\ge\frac{2688}{5}>480=T_6,\\[4pt]
N=7:&\displaystyle\Lambda_1(B_1)\ge504>
\frac{4851}{19}=T_7.
\end{array}
\]
Thus every ball satisfies \eqref{eq:domaincondition} in dimensions
$5,6,7$, and \textup{(i)} gives (H). For $N\ge8$, (H) is not
required. Consequently Theorem~\ref{thm:main} applies on every ball
for all $N\ge5$, with the second solution existing on the set of
masses specified there.

\medskip
\noindent\textup{(iii)} We exhibit non-star-shaped domains satisfying
\eqref{eq:domaincondition}. Suppose first that
\[
B_R(x_0)\subset\Omega\subset B_{R'}(x_0),\qquad 0<R<R'.
\]
Zero extension gives $H_0^2(\Omega)\subset H_0^2(B_{R'}(x_0))$,
so the Rayleigh characterization implies
\[
\Lambda_1(\Omega)\ge\Lambda_1(B_{R'}(x_0))
=(R')^{-4}\Lambda_1(B_1).
\]
Since $R_\Omega\ge R$, we obtain
\[
\Lambda_1(\Omega)R_\Omega^4\ge
\begin{cases}
500(R/R')^4,&N=5,\\
(2688/5)(R/R')^4,&N=6,\\
504(R/R')^4,&N=7.
\end{cases}
\]
It is therefore sufficient to choose $R,R'$ with
\begin{equation}\label{eq:annularcondition}
1<\frac{R'}{R}<
\begin{cases}
(500/T_5)^{1/4},&N=5,\\
((2688/5)/T_6)^{1/4},&N=6,\\
(504/T_7)^{1/4},&N=7.
\end{cases}
\end{equation}
Each upper bound is greater than one, so such radii exist.
Choose a closed ball $\overline{B_r(y)}$ contained in
$B_{R'}(x_0)\setminus\overline{B_R(x_0)}$, and set
\[
\Omega=B_{R'}(x_0)\setminus\overline{B_r(y)}.
\]
This domain is smooth, bounded and connected, and satisfies the
required inclusions. Hence \eqref{eq:annularcondition} implies
\eqref{eq:domaincondition}.
It is not star-shaped with respect to any point. Indeed, for any
$z\in\Omega$, the ray from $z$ through $y$ crosses the removed ball
and then meets $\Omega$ again before leaving $B_{R'}(x_0)$, since
$\overline{B_r(y)}\Subset B_{R'}(x_0)$. Choosing $w$ in that latter
portion gives $z,w\in\Omega$ but $[z,w]\not\subset\Omega$.
This proves \textup{(iii)}.
\end{proof}

\end{document}